\documentclass[12pt]{article}
\begin{document}

\author{Md Fazlul Hoque\\School of Mathematics and Physics\\The University of Queensland, Brisbane, Australia\\E-mail: m.hoque@uq.edu.au
\\\\A.C. Paul\\Department of Mathematics\\Rajshahi University, Bangladesh\\E-mail: acpaulrubd\_math@yahoo.com}

\title{Centralizers on Prime and Semiprime Gamma Rings}
\maketitle

\newtheorem{th1}{Theorem}[section]
\newtheorem{df}{Definition}
\newtheorem{lem}{Lemma}[section]
\newtheorem{cor}{Corollary}[section]
\newtheorem{pro}{Proposition}[section]

\begin{abstract}Let $M$ be a noncommutative 2-torsion free semiprime $\Gamma$-ring satisfying a certain assumption and let $S$ and $T$ be left centralizers on $M$. We prove the following results:
\\(i) If $[S(x),T(x)]_{\alpha }\beta S(x)+S(x)\beta [S(x),T(x)]_{\alpha }$=$0$ holds for all $x\in M$ and $\alpha ,\beta \in \Gamma $, then $[S(x),T(x)]_{\alpha }$=$0$.
\\(ii) If $S\neq 0 (T\neq 0)$, then there exists $\lambda \in C$,(the extended centroid of $M$) such that $T$=$\lambda \alpha S(S=\lambda \alpha T)$ for all $\alpha \in \Gamma $.
\\(iii) Suppose that $[[S(x),T(x)]_{\alpha },S(x)]_{\beta }$=$0$ holds for all $x\in M$ and $\alpha ,\beta \in\Gamma $. Then $[S(x),T(x)]_{\alpha }$=$0$ for all $x\in M$ and $\alpha \in\Gamma $.
\\(iv) If $M$ is a prime $\Gamma $-ring satisfying a certain assumption and $S\neq 0(T\neq 0)$, then there exists $\lambda \in C$, the extended centroid, such that $T$=$\lambda \alpha S(S=\lambda \alpha T)$. 
\end{abstract}
\section{Introduction}The notion of a $\Gamma $-ring was first introduced as an extensive generalization of the concept of a classical ring. From its first appearance, the extensions and generalizations of various important results in the theory of classical rings to the theory of $\Gamma $-rings have been attracted a wider attentions as an emerging area of research to the modern algebraists to enrich the world of algebra. All over the world, many prominent mathematicians have worked out on this interesting area of research to determine many basic properties of $\Gamma $-rings and have  executed more productive and creative results of $\Gamma $-rings in the last few decades. We begin with the definition.
\par
%
%
Let $M$ and $\Gamma $  be additive abelian groups. If there exists an additive mapping $(x,\alpha ,y)\rightarrow  x\alpha y$ of $M\times\Gamma \times M\rightarrow M$, which satisfies the conditions $(x\alpha y)\beta z$=$x\alpha (y\beta z)$ for all $x,y,z\in M$ and $\alpha ,\beta \in \Gamma$, then $M$ is called a $\Gamma$-ring.
Every ring $M$ is a $\Gamma $-ring with $M$=$\Gamma$. However a $\Gamma $-ring need not be a ring. Gamma rings, more general than rings[8]. Barnes[1] weakened slightly the conditions in the definition of $\Gamma $-ring in the sense of Nobusawa. 
Let $M$ be a $\Gamma $-ring. Then an additive subgroup $U$ of $M$ is called a left (right) ideal of $M$ if $M\Gamma U\subset U$($U\Gamma M \subset U$). If $U$ is both a left and a right ideal, then we say $U$ is an ideal of $M$.
Suppose again that $M$ is a $\Gamma $-ring. Then $M$ is said to be a 2-torsion free if $2x$=$0$ implies $x$=$0$ for all $x\in M$.
An ideal $P_{1}$ of a $\Gamma $-ring $M$ is said to be prime if for any ideals $A$ and $B$ of $M$, $A\Gamma B\subseteq P_{1}$ implies $A\subseteq P_{1}$ or $B\subseteq P_{1}$. An ideal $P_{2}$ of a $\Gamma $-ring $M$ is said to be semiprime if for any ideal $U$ of $M$, $U\Gamma U\subseteq P_{2}$ implies $U\subseteq P_{2}$. A $\Gamma$-ring $M$ is said to be prime if $a\Gamma M\Gamma b$=$(0)$ with $a,b\in M$, implies $a$=$0$ or $b$=$0$ and semiprime if $a\Gamma M\Gamma a$=$(0)$ with $a\in M$ implies $a$=$0$. Furthermore, $M$ is said to be commutative $\Gamma $-ring if $x\alpha y$=$y\alpha x$ for all $x, y\in M$ and $\alpha \in\Gamma $. 
Moreover, the set $Z(M)$ =$\{x\in M:x\alpha y=y\alpha x $ for all $ \alpha \in \Gamma, y\in M\}$ is called the centre of the $\Gamma $-ring $M$.
For the definitions of the quotent $\Gamma $-ring, the extended centroid and the central closure, we refer to [7].
\\If $M$ is a $\Gamma $-ring, then $[x,y]_{\alpha }$=$x\alpha y-y\alpha x$ is known as the commutator of $x$ and $y$ with respect to $\alpha $, where $x,y\in M$ and $\alpha \in\Gamma $. We make the basic commutator identities:
\\  $[x\alpha y,z]_{\beta }$=$[x,z]_{\beta }\alpha y+x[\alpha ,\beta ]_{z}y+x\alpha [y,z]_{\beta }$ \\and  $[x,y\alpha z]_{\beta }$=$[x,y]_{\beta }\alpha z+y[\alpha ,\beta ]_{x}z+y\alpha [x,z]_{\beta }$ for all $x,y.z\in M$ and $\alpha ,\beta \in\Gamma $.
We consider the following assumption:\\$(A)$.................$x\alpha y\beta z$=$x\beta y\alpha z$ for all $x,y,z\in M$ and $\alpha ,\beta \in\Gamma $.\\ According to the assumption $(A)$, the above two identites reduce to \\ $[x\alpha y,z]_{\beta }$=$[x,z]_{\beta }\alpha y+x\alpha [y,z]_{\beta }$ \\and  $[x,y\alpha z]_{\beta }$=$[x,y]_{\beta }\alpha z+y\alpha [x,z]_{\beta }$, which we extensively used.
\\An additive mapping $T : M\rightarrow M$ is a left(right) centralizer if $T(x\alpha y)$=$T(x)\alpha y$ $(T(x\alpha y)=x\alpha T(y))$ holds for all $x,y\in M$ and $\alpha\in \Gamma$. A centralizer is an additive mapping which is both a left and a right centralizer. For any fixed $a\in M$ and $\alpha \in\Gamma$, the mapping $T(x)=a\alpha x$ is a left centralizer and $T(x)=x\alpha a$ is a right centralizer. We shall restrict our attention on left centralizer since all results represented in this paper are true also for right centralizers because of left and right symmetry.
\\ Barnes [1], Lue [6] and Kyuno[5] studied the structure of $\Gamma$-rings and obtained various generalizations of corresponding parts in ring theory.
\\ Borut Zalar[14] worked on centralizers of semiprime rings and proved that Jordan centralizers and centralizers of these rings coincide. Joso Vukman[11, 12, 13] developed some remarkable results using centralizers on prime and semiprime rings.
\\In [2], Hoque and Paul have proved that every Jordan centralizer of a 2-torsion free semiprime $\Gamma $-ring satisfying a certain assumption is a centralizer. 
Also, they proved in [3],  if $T$ is an additive mapping on a 2-torsion free semiprime $\Gamma $-ring $M$ with a certain assumption such that $T(x\alpha y\beta x)=x\alpha T(y)\beta x$ for all $x, y\in M$ and $\alpha ,\beta \in\Gamma $, then $T$ is a centralizer
 and in [4], if $2T(x\alpha b\beta a)=T(x)\alpha y\beta x+x\alpha y\beta T(x)$ holds for all $x, y\in M$ and $\alpha ,\beta \in\Gamma $, then $T$ is also a centralizer.
\\ In this paper, we generalize some results of Joso Vukman[11] in Gamma rings.
  
\section{Centralizers of Prime and Semiprime $\Gamma $-rings.}
  To prove our main results, we need the following lammas:
\vskip.2cm
\noindent
\begin{lem} Suppose that the elements  $a_{i}$, $b_{i}$ in the central closure of a prime $\Gamma $-ring $M$ satisfy $\Sigma a_{i}\alpha _{i}x\beta _{i}b_{i}=0$ for all $x\in M$ and $\alpha _{i},\beta _{i}\in\Gamma $. If $b_{i}\neq 0$ for some $i$, then $a_{i}$'s are $C$-dependent, where $C$ is the extended centroid. 
\end{lem}
\vskip.2cm
\noindent
{\bf Proof.} Let $M$ be a prime $\Gamma $-ring and let $C_{\Gamma }=C$ be the extended centroid of $M$. If $a_{i}$ and $b_{i}$ are non-zero elements of $M$ such that $\Sigma a_{i}\alpha _{i}x\beta _{i}b_{i}=0$ for all $x\in M$ and $\alpha _{i},\beta _{i}\in\Gamma $, then $a_{i}$'s (also $b_{i}$'s) are linearly dependent over $C$. Moreover, if $a\alpha x\beta b=b\alpha x\beta a$ for all $x\in M$ and $\alpha ,\beta \in\Gamma $, where $a(\neq 0)$, $b\in M$ are fixed, then there exists $\lambda \in C$ such that $a=\lambda \alpha b$ for $\alpha \in\Gamma $. Clearly, the lemma is proved.   
\vskip.2cm
\noindent
\begin{lem} Suppose that $M$ is a noncommutative prime $\Gamma$-ring satisfying the assumption $(A)$ and $T:M\rightarrow M$ is a left centralizer. If $T(x)\in Z(M)$ for all $x\in M$, then  $T=0$.
\end{lem}
\vskip.2cm
\noindent
{\bf Proof.} Since $T$ is a left centralizer on $M$, we have $T(x\alpha y)=T(x)\alpha y$ holds for all $x,y\in M$ and $\alpha \in\Gamma $ and hence $[T(x),y]_{\alpha }=0$ for all $x,y\in M$ and $\alpha \in\Gamma $.
Putting $x=x\beta z$ in the above relation, we have \begin{eqnarray*}&0&=[T(x\beta z),y]_{\alpha }\\&&=[T(x)\beta z,y]_{\alpha }\\&&=[T(x),y]_{\alpha }\beta z+T(x)\beta [z,y]_{\alpha }\\&&=T(x)\beta [z,y]_{\alpha }\end{eqnarray*}
Hence $T(x)\beta [z,y]_{\alpha }=0$, which gives $T(x)\beta w \gamma [z,y]_{\alpha }=0$ for all $x,y.z,w\in M$ and $\alpha ,\beta ,\gamma \in\Gamma $, whence it follows that $T=0$, otherwise $M=0$.
\vskip.2cm
\noindent
\begin{lem} Suppose that $M$ is a noncommutative prime $\Gamma$-ring satisfying the assumption $(A)$ and $S, T:M\rightarrow M$ are left centralizers. If $[S(x),T(x)]_{\alpha }=0$ holds for all $x\in M$ and $\alpha \in\Gamma $ and $T\neq 0$, then there exists $\lambda \in C$ such that $S=\lambda \alpha T$, where $C$ is the extended centroid. 
\end{lem}
\vskip.2cm
\noindent
{\bf Proof.} First, we put $x=x+y$ in the relation $[S(x),T(x)]_{\alpha }=0$ and linearize, we have \begin{eqnarray}&[S(x),T(y)]_{\alpha }+[S(y),T(x)]_{\alpha }&=0\end{eqnarray}. 
Replace $y$ by $y\beta z$ in (1), we have \begin{eqnarray*}&0&=[S(x),T(y)\beta z]_{\alpha }+[S(y)\beta z,T(x)]_{\alpha }\\&&=[S(x),T(y)]_{\alpha }\beta z+T(y)\beta [S(x),z]_{\alpha }+[S(y),T(x)]_{\alpha }\beta z+S(y)\beta [z,T(x)]_{\alpha }\\&&=T(y)\beta [S(x),z]_{\alpha }+S(y)\beta [z,T(x)]_{\alpha }\end{eqnarray*} Thus we have \begin{eqnarray*}&T(y)\beta [S(x),z]_{\alpha }+S(y)\beta [z,T(x)]_{\alpha }&=0\end{eqnarray*}
Putting $y=y\gamma w$ in the above relation, we obtain, \begin{eqnarray}&T(y)\gamma w\beta [S(x),z]_{\alpha }+S(y)\gamma w\beta [z,T(x)]_{\alpha }&=0\end{eqnarray} Our assumption $T\neq 0$ follows from Lemma-2.2 that there exist $x,z\in M$ and $\alpha \in \Gamma $ such that $[T(x),z]_{\alpha }\neq 0$. Now, the relation (2) and Lemma-2.1 imply that $S(y)=\lambda (y)\alpha T(y)$, where $ \lambda (y)$ is from $C$.
If we put $S(y)=\lambda (y)\alpha T(y)$ and $S(x)=\lambda (x)\alpha T(x)$ in the relation (2), we obtain $(\lambda (x)-\lambda (y))\alpha T(y)\gamma w\beta [T(x),z]_{\alpha }=0$ for all pairs $y,w\in M$, whence it follows $(\lambda (x)-\lambda (y))\alpha T(y)=0$, since $[T(x),z]_{\alpha }\neq 0$. Thus we have $\lambda (x)\alpha T(y)=\lambda (y)\alpha T(y)$ which completes the proof of the lemma.
\vskip.2cm
\noindent
\begin{th1} Suppose that $M$ is a  2-torsion free noncommutative semiprime $\Gamma $-ring satisfying the assumption $(A)$ and $S$, $T$ are left centralizers on $M$. If $[S(x),T(x)]_{\alpha }\beta S(x)+S(x)\beta [S(x),T(x)]_{\alpha }=0$ holds for all $x\in M$ and $\alpha ,\beta \in\Gamma $. Then $[S(x),T(x)]_{\alpha }=0$ for all $x\in M$ and $\alpha \in\Gamma $. 
Also, if $M$ is prime $\Gamma $-ring satisfying the assumption $(A)$ and $S\neq 0(T\neq 0)$, then there esixts $\lambda \in C$,(the extended centroid of $M$) such that $T=\lambda \alpha S (S=\lambda \alpha T)$.
\end{th1}
\vskip.2cm
\noindent
{\bf Proof.}  By the hypothesis, we have \begin{eqnarray}&[S(x),T(x)]_{\alpha }\beta S(x)+S(x)\beta [S(x),T(x)]_{\alpha }&=0\end{eqnarray} The lineariztion of the above relation, we have \begin{eqnarray}0&=&[S(x),T(x)]_{\alpha }\beta S(y)+S(y)\beta [S(x),T(x)]_{\alpha }\nonumber\\&&+[S(x),T(y)]_{\alpha }\beta S(x)+S(x)\beta [S(x),T(y)]_{\alpha }\nonumber\\&&+[S(y),T(x)]_{\alpha }\beta S(x)+S(x)\beta [S(y),T(x)]_{\alpha }\nonumber\\&&+[S(y),T(y)]_{\alpha }\beta S(x)+S(x)\beta [S(y),T(y)]_{\alpha }\nonumber\\&&+[S(y),T(x)]_{\alpha }\beta S(y)+S(y)\beta [S(y),T(x)]_{\alpha }\nonumber\\&&+[S(x),T(y)]_{\alpha }\beta S(y)+S(y)\beta [S(x),T(y)]_{\alpha }\end{eqnarray}
Replacing $-x$ for $x$ in the above relation, we have \begin{eqnarray}0&=&[S(x),T(x)]_{\alpha }\beta S(y)+S(y)\beta [S(x),T(x)]_{\alpha }\nonumber\\&&+[S(x),T(y)]_{\alpha }\beta S(x)+S(x)\beta [S(x),T(y)]_{\alpha }\nonumber\\&&+[S(y),T(x)]_{\alpha }\beta S(x)+S(x)\beta [S(y),T(x)]_{\alpha }\nonumber\\&&-[S(y),T(y)]_{\alpha }\beta S(x)-S(x)\beta [S(y),T(y)]_{\alpha }\nonumber\\&&-[S(y),T(x)]_{\alpha }\beta S(y)-S(y)\beta [S(y),T(x)]_{\alpha }\nonumber\\&&-[S(x),T(y)]_{\alpha }\beta S(y)-S(y)\beta [S(x),T(y)]_{\alpha }\end{eqnarray}
Adding (4) and (5), we have \begin{eqnarray*}0&=&2[S(x),T(x)]_{\alpha }\beta S(y)+2S(y)\beta [S(x),T(x)]_{\alpha }\nonumber\\&&+2[S(x),T(y)]_{\alpha }\beta S(x)+2S(x)\beta [S(x),T(y)]_{\alpha }\nonumber\\&&+2[S(y),T(x)]_{\alpha }\beta S(x)+2S(x)\beta [S(y),T(x)]_{\alpha }\end{eqnarray*} Hence by 2-torsion freeness of $M$, it follows that \begin{eqnarray}0&=&[S(x),T(x)]_{\alpha }\beta S(y)+S(y)\beta [S(x),T(x)]_{\alpha }\nonumber\\&&+[S(x),T(y)]_{\alpha }\beta S(x)+S(x)\beta [S(x),T(y)]_{\alpha }\nonumber\\&&+[S(y),T(x)]_{\alpha }\beta S(x)+S(x)\beta [S(y),T(x)]_{\alpha }\end{eqnarray}
Replacing $y$ by $x\gamma y$ in the above relation , we have \begin{eqnarray*}0&=&[S(x),T(x)]_{\alpha }\beta S(x)\gamma y+S(x)\gamma y\beta [S(x),T(x)]_{\alpha }\nonumber\\&&+[S(x),T(x)\gamma y]_{\alpha }\beta S(x)+S(x)\beta [S(x),T(x)\gamma y]_{\alpha }\nonumber\\&&+[S(x)\gamma y,T(x)]_{\alpha }\beta S(x)+S(x)\beta [S(x)\gamma y,T(x)]_{\alpha }\nonumber\\&=&[S(x),T(x)]_{\alpha }\beta S(x)\gamma y+S(x)\gamma y\beta [S(x),T(x)]_{\alpha }\nonumber\\&&+[S(x),T(x)]_{\alpha }\gamma y\beta S(x)+T(x)\gamma [S(x),y]_{\alpha }\beta S(x)\\&&+S(x)\beta [S(x),T(x)]_{\alpha }\gamma y+S(x)\beta T(x)\gamma [S(x),y]_{\alpha }\nonumber\\&&+[S(x),T(x)]_{\alpha }\gamma y\beta S(x)+S(x)\gamma [y,T(x)]_{\alpha }\beta S(x)\\&&+S(x)\beta [S(x),T(x)]_{\alpha }\gamma y+S(x)\beta S(x)\gamma [y,T(x)]_{\alpha }\end{eqnarray*}
According to (6), the above relation reduces to \begin{eqnarray}0&=&S(x)\gamma y\beta [S(x),T(x)]_{\alpha }+2[S(x),T(x)]_{\alpha }\gamma y\beta S(x)\nonumber\\&&+T(x)\gamma [S(x),y]_{\alpha }\beta S(x)+S(x)\beta T(x)\gamma [S(x),y]_{\alpha }\nonumber\\&&+S(x)\gamma [y,T(x)]_{\alpha }\beta S(x)+S(x)\beta [S(x),T(x)]_{\alpha }\gamma y\nonumber\\&&+S(x)\beta S(x)\gamma [y,T(x)]_{\alpha }\end{eqnarray}
Putting $y=y\delta S(x)$ in (7), we obtain \begin{eqnarray*}0&=&S(x)\gamma y\delta S(x)\beta [S(x),T(x)]_{\alpha }+2[S(x),T(x)]_{\alpha }\gamma y\delta S(x)\beta S(x)\nonumber\\&&+T(x)\gamma [S(x),y]_{\alpha }\beta S(x)\delta S(x)+S(x)\beta T(x)\gamma [S(x),y]_{\alpha }\delta S(x)\nonumber\\&&+S(x)\gamma [y,T(x)]_{\alpha }\delta S(x)\beta S(x)+S(x)\gamma y\delta [S(x),T(x)]_{\alpha }\beta S(x)\nonumber\\&&+S(x)\beta [S(x),T(x)]_{\alpha }\gamma y\delta S(x)+S(x)\beta S(x)\gamma y\delta [S(x),T(x)]_{\alpha }\nonumber\\&&+S(x)\beta S(x)\gamma [y,T(x)]_{\alpha }\delta S(x)\end{eqnarray*}
which gives according to (7) to \begin{eqnarray}0&=&S(x)\gamma y\delta S(x)\beta [S(x),T(x)]_{\alpha }+S(x)\beta S(x)\gamma y\delta [S(x),T(x)]_{\alpha }\end{eqnarray} Putting $y=T(x)\omega y$ in (8), we have \begin{eqnarray}0&=&S(x)\gamma T(x)\omega y\delta S(x)\beta [S(x),T(x)]_{\alpha }\nonumber\\&&+S(x)\beta S(x)\gamma T(x)\omega y\delta [S(x),T(x)]_{\alpha }\end{eqnarray} Also left multiplication of (8) by $T(x)\omega $ gives \begin{eqnarray}0&=&T(x)\omega S(x)\gamma y\delta S(x)\beta [S(x),T(x)]_{\alpha }\nonumber\\&&+T(x)\omega S(x)\beta S(x)\gamma y\delta [S(x),T(x)]_{\alpha }\end{eqnarray} 
From (9) and (10), we obtain, \begin{eqnarray*}0&=&[S(x),T(x)]_{\gamma }\omega y\delta S(x)\beta [S(x),T(x)]_{\alpha }+[S(x)\beta S(x),T(x)]_{\gamma }\omega y\delta [S(x),T(x)]_{\alpha }\nonumber\\&=&[S(x),T(x)]_{\gamma }\omega y\delta S(x)\beta [S(x),T(x)]_{\alpha }\nonumber\\&&+([S(x),T(x)]_{\gamma }\beta S(x)+S(x)\beta [S(x),T(x)]_{\gamma })\omega y\delta [S(x),T(x)]_{\alpha }\nonumber\\&=&[S(x),T(x)]_{\gamma }\omega y\delta S(x)\beta [S(x),T(x)]_{\alpha }\end{eqnarray*}   
Thus we have \begin{eqnarray*}0&=&[S(x),T(x)]_{\gamma }\omega y\delta S(x)\beta [S(x),T(x)]_{\alpha }\end{eqnarray*} Left multiplication of the above relation by $S(x)\beta $ gives \begin{eqnarray}0&=&S(x)\beta [S(x),T(x)]_{\gamma }\omega y\delta S(x)\beta [S(x),T(x)]_{\alpha }\end{eqnarray} for all $x,y\in M$ and $\alpha ,\beta ,\gamma ,\delta ,\omega \in\Gamma $.
Hence from (11), it follows \begin{eqnarray}S(x)\beta [S(x),T(x)]_{\alpha }&=&0\end{eqnarray} From (3) and (12), we have also \begin{eqnarray}[S(x),T(x)]_{\alpha }\beta S(x)&=&0\end{eqnarray} From (12) one obtains the relation \begin{eqnarray}0&=&S(y)\beta [S(x),T(x)]_{\alpha }+S(x)\beta [S(y),T(x)]_{\alpha }\nonumber\\&&+S(x)\beta [S(x),T(y)]_{\alpha }\end{eqnarray}(see the proof of (6)). 
Putting $y=x\gamma y$ in (14), we have \begin{eqnarray*}0&=&S(x)\gamma y\beta [S(x),T(x)]_{\alpha }+S(x)\beta [S(x)\gamma y,T(x)]_{\alpha }\nonumber\\&&+S(x)\beta [S(x),T(x)\gamma y]_{\alpha }\\&=&S(x)\gamma y\beta [S(x),T(x)]_{\alpha }+S(x)\beta [S(x),T(x)]_{\alpha }\gamma y\nonumber\\&&+S(x)\beta S(x)\gamma [y,T(x)]_{\alpha }+S(x)\beta [S(x),T(x)]_{\alpha }\gamma y+S(x)\beta T(x)\gamma [S(x),y]_{\alpha }\\&=&S(x)\gamma y\beta [S(x),T(x)]_{\alpha }+S(x)\beta S(x)\gamma [y,T(x)]_{\alpha }\\&&+S(x)\beta T(x)\gamma [S(x),y]_{\alpha }\end{eqnarray*}
Thus we have the above relation  \begin{eqnarray*}0&=&S(x)\gamma y\beta [S(x),T(x)]_{\alpha }+S(x)\beta S(x)\gamma [y,T(x)]_{\alpha }\\&&+S(x)\beta T(x)\gamma [S(x),y]_{\alpha }\end{eqnarray*} which can be written in the form \begin{eqnarray*}&0=&S(x)\gamma y\beta [S(x),T(x)]_{\alpha }+S(x)\beta S(x)\gamma y\alpha T(x)\\&&-S(x)\beta T(x)\gamma y\alpha S(x)+S(x)\beta [T(x),S(x)]_{\gamma }\alpha y\end{eqnarray*} whence it follows \begin{eqnarray}&0=&S(x)\gamma y\beta [S(x),T(x)]_{\alpha }+S(x)\beta S(x)\gamma y\alpha T(x)\nonumber\\&&-S(x)\beta T(x)\gamma y\alpha S(x)\end{eqnarray} according to (12). 
Taking $T(x)\delta $ of (15) on the left side, we have \begin{eqnarray}&0=&T(x)\delta S(x)\gamma y\beta [S(x),T(x)]_{\alpha }+T(x)\delta S(x)\beta S(x)\gamma y\alpha T(x)\nonumber\\&&-T(x)\delta S(x)\beta T(x)\gamma y\alpha S(x)\end{eqnarray} Putting $y=T(x)\delta y$ in (15) gives \begin{eqnarray}&0=&S(x)\gamma T(x)\delta y\beta [S(x),T(x)]_{\alpha }+S(x)\beta S(x)\gamma T(x)\delta y\alpha T(x)\nonumber\\&&-S(x)\beta T(x)\gamma T(x)\delta y\alpha S(x)\end{eqnarray}
From (16) and (17), we have \begin{eqnarray*}&0=&[S(x),T(x)]_{\gamma }\delta y\beta [S(x),T(x)]_{\alpha }+[S(x)\beta S(x),T(x)]_{\gamma }\delta y\alpha T(x)\nonumber\\&&+[T(x),S(x)]_{\delta }\beta T(x)\gamma y\alpha S(x)\\&=&[S(x),T(x)]_{\gamma }\delta y\beta [S(x),T(x)]_{\alpha }+([S(x),T(x)]_{\gamma }\beta S(x)\\&&+S(x)\beta [S(x),T(x)]_{\gamma })\delta y\alpha T(x)+[T(x),S(x)]_{\delta }\beta T(x)\gamma y\alpha S(x)\end{eqnarray*} which reduces to \begin{eqnarray}0&=&[S(x),T(x)]_{\gamma }\delta y\beta [S(x),T(x)]_{\alpha }+[T(x),S(x)]_{\delta }\beta T(x)\gamma y\alpha S(x)\end{eqnarray}    
The substitution $y\omega S(x)\rho z$ for $y$ in (18) gives \begin{eqnarray}0&=&[S(x),T(x)]_{\gamma }\delta y\omega S(x)\rho z\beta [S(x),T(x)]_{\alpha }\nonumber\\&&+[T(x),S(x)]_{\delta }\beta T(x)\gamma y\omega S(x)\rho z\alpha S(x)\end{eqnarray} Again, right multiplication of (18) by $\omega z\rho S(x)$, we have  \begin{eqnarray}0&=&[S(x),T(x)]_{\gamma }\delta y\beta [S(x),T(x)]_{\alpha }\omega z\rho S(x)\nonumber\\&&+[T(x),S(x)]_{\delta }\beta T(x)\gamma y\alpha S(x)\omega z\rho S(x)\end{eqnarray} 
From (19) and (20), we obtain \begin{eqnarray}0&=&[S(x),T(x)]_{\gamma }\delta y\beta A(x,z)\end{eqnarray} where $A(x,z)=[S(x),T(x)]_{\alpha }\omega z\rho S(x)-S(x)\rho z\omega [S(x).T(x)]_{\alpha }$. Replacing $y$ by $z\rho S(x)\omega y$ in (21) gives \begin{eqnarray}0&=&[S(x),T(x)]_{\gamma }\delta z\rho S(x)\omega y\beta A(x,z)\end{eqnarray}  Left multiplication of (21) by $S(x)\rho z\omega $ gives \begin{eqnarray}0&=&S(x)\rho z\omega [S(x),T(x)]_{\gamma }\delta y\beta A(x,z)\end{eqnarray} Combining (22) and (23), we arrive at \begin{eqnarray*}0&=&A(x,z)\delta y\beta A(x,z)\end{eqnarray*} for all $x,y,z\in M$ and $\delta ,\beta \in \Gamma $.
Hence by semiprimeness of $M$, it follows $A(x,z)=0$ and hence \begin{eqnarray}[S(x),T(x)]_{\alpha }\omega z\rho S(x)&=&S(x)\rho z\omega [S(x),T(x)]_{\alpha }\end{eqnarray} The substitution of $z$ by $T(x)\gamma y$ in (24) gives  \begin{eqnarray}[S(x),T(x)]_{\alpha }\omega T(x)\gamma y\rho S(x)&=&S(x)\rho T(x)\gamma y\omega [S(x),T(x)]_{\alpha }\end{eqnarray} The relation (25) makes it possible to replace in (18), $[S(x),T(x)]_{\delta }\beta T(x)\gamma y\alpha S(x)$ by $S(x)\alpha T(x)\gamma y\beta [S(x),T(x)]_{\delta }$. Thus we have  \begin{eqnarray*}0&=&[S(x),T(x)]_{\gamma }\delta y\beta [S(x),T(x)]_{\alpha }-S(x)\alpha T(x)\gamma y\beta [S(x),T(x)]_{\delta }\end{eqnarray*} which reduces to \begin{eqnarray}0&=&T(x)\gamma S(x)\delta y\beta [S(x),T(x)]_{\alpha }\end{eqnarray}  
Putting  $y=T(x)\omega y$ in (26), we have \begin{eqnarray}0&=&T(x)\gamma S(x)\delta T(x)\omega y\beta [S(x),T(x)]_{\alpha }\end{eqnarray} Multiplying (26) from the left side by $T(x)\omega $, we have \begin{eqnarray}0&=&T(x)\omega T(x)\gamma S(x)\delta y\beta [S(x),T(x)]_{\alpha }\end{eqnarray} Subtrating (28) from (27), we have \begin{eqnarray*}0&=&T(x)\omega [S(x),T(x)]_{\gamma }\delta y\beta [S(x),T(x)]_{\alpha }\end{eqnarray*} which gives putting $y=y\omega T(x)$,\begin{eqnarray*}0&=&T(x)\omega [S(x),T(x)]_{\gamma }\delta y\omega T(x)\beta [S(x),T(x)]_{\alpha }\\&=&T(x)\omega [S(x),T(x)]_{\gamma }\delta y\beta T(x)\omega  [S(x),T(x)]_{\alpha }\end{eqnarray*} whence it follows \begin{eqnarray}0&=&T(x)\omega [S(x),T(x)]_{\alpha }\end{eqnarray}             
The substitution $y=y\beta T(x)$ in (25) gives because of (29)\begin{eqnarray}0&=&[S(x),T(x)]_{\alpha }\omega y\beta T(x)\rho S(x)\end{eqnarray} From (13), we obtain the relation \begin{eqnarray*}0&=&[S(x),T(x)]_{\alpha }\beta S(y)+[S(x),T(y)]_{\alpha }\beta S(x)+[S(y),T(x)]_{\alpha }\beta S(x)\end{eqnarray*} (see the proof of (6)). Putting in the above relation $y=x\gamma y$, we have \begin{eqnarray*}0&=&[S(x),T(x)]_{\alpha }\beta S(x)\gamma y+[S(x),T(x)\gamma y]_{\alpha }\beta S(x)+[S(x)\gamma y,T(x)]_{\alpha }\beta S(x)\\&=&T(x)\gamma [S(x),y]_{\alpha }\beta S(x)+[S(x),T(x)]_{\alpha }\gamma y\beta S(x)\\&&+S(x)\gamma [y,T(x)]_{\alpha }\beta S(x)+[S(x),T(x)]_{\alpha }\gamma y\beta S(x)\end{eqnarray*} 
Thus we have   \begin{eqnarray*}&0=&2[S(x),T(x)]_{\alpha }\gamma y\beta S(x)+T(x)\gamma [S(x),y]_{\alpha }\beta S(x)+S(x)\gamma [y,T(x)]_{\alpha }\beta S(x)\end{eqnarray*} which can be written after some calculation in the form \begin{eqnarray}&0=&[S(x),T(x)]_{\alpha }\gamma y\beta S(x)+S(x)\gamma y\alpha T(x)\beta S(x)\nonumber\\&&-T(x)\gamma y\alpha S(x)\beta S(x)\end{eqnarray} The relation (24) makes it possible to replace in (31), $[S(x),T(x)]_{\alpha }\gamma y\beta S(x)$ by $S(x)\beta y\gamma [S(x),T(x)]_{\alpha }$. Thus we have  \begin{eqnarray*}&0=&S(x)\beta y\gamma [S(x),T(x)]_{\alpha }+S(x)\gamma y\alpha T(x)\beta S(x)\nonumber\\&&-T(x)\gamma y\alpha S(x)\beta S(x)\\&=&S(x)\beta y\gamma S(x)\alpha T(x)-T(x)\gamma y\alpha S(x)\beta S(x)\end{eqnarray*}
Therefore, we have \begin{eqnarray}S(x)\beta y\gamma S(x)\alpha T(x)&=&T(x)\gamma y\alpha S(x)\beta S(x)\end{eqnarray} Putting in the above relation $y=T(x)\omega y$, we have \begin{eqnarray}S(x)\beta T(x)\omega y\gamma S(x)\alpha T(x)&=&T(x)\gamma T(x)\omega y\alpha S(x)\beta S(x)\end{eqnarray} Left multiplication of (32) by $T(x)\omega $ gives \begin{eqnarray}T(x)\omega S(x)\beta y\gamma S(x)\alpha T(x)&=&T(x)\omega T(x)\gamma y\alpha S(x)\beta S(x)\end{eqnarray} Combining (33) and (34), we have \begin{eqnarray*}0&=&[S(x),T(x)]_{\beta }\omega y\gamma S(x)\alpha T(x)\end{eqnarray*} which gives together with (30), \begin{eqnarray*}0&=&[S(x),T(x)]_{\alpha }\omega y\beta [S(x),T(x)]_{\alpha }\end{eqnarray*} Hence by semiprimeness of $M$, we have \begin{eqnarray}[S(x),T(x)]_{\alpha }&=&0\end{eqnarray}

If $M$ is a prime $\Gamma $-ring, then the relation (35) and Lemma-2.3 complete the proof of the theorem.      
          
\vskip.2cm
\noindent
\begin{th1} Suppose that $M$ is a 2-torsion free noncommutative semiprime $\Gamma $-ring satisfying the assumption $(A)$ and $S$, $T$ are left centralizers on $M$. If $[[S(x),T(x)]_{\alpha }, S(x)]_{\beta }=0$ holds for all $x\in M$ and $\alpha ,\beta \in\Gamma $. Then $[S(x),T(x)]_{\alpha }=0$ for all $x\in M$ and $\alpha \in\Gamma $. 
Moreover, if $M$ is prime $\Gamma $-ring satisfying the assumption $(A)$ and $S\neq 0(T\neq 0)$, then there esixts $\lambda \in C$,(the extended centroid of $M$) such that $T=\lambda \alpha S (S=\lambda \alpha T)$.
\end{th1}
\vskip.2cm
\noindent
{\bf Proof.} By the assumption \begin{eqnarray}[[S(x),T(x)]_{\alpha },S(x)]_{\beta }&=&0\end{eqnarray} The linearization of (36) gives \begin{eqnarray}0&=&[[S(x),T(x)]_{\alpha },S(y)]_{\beta }+[[S(x),T(y)]_{\alpha },S(x)]_{\beta }\nonumber\\&&+[[S(y),T(x)]_{\alpha },S(x)]_{\beta }\end{eqnarray} Putting $y=x\gamma y$ in (37), we have \begin{eqnarray*}0&=&[[S(x),T(x)]_{\alpha },S(x)\gamma y]_{\beta }+[[S(x),T(x)\gamma y]_{\alpha },S(x)]_{\beta }\nonumber\\&&+[[S(x)\gamma y,T(x)]_{\alpha },S(x)]_{\beta }\\&=&[[S(x),T(x)]_{\alpha },S(x)]_{\beta }\gamma y+S(x)\gamma [[S(x),T(x)]_{\alpha },y]_{\beta }\nonumber\\&&+[[S(x),T(x)]_{\alpha }\gamma y+T(x)\gamma [S(x),y]_{\alpha },S(x)]_{\beta }\\&&+[[S(x),T(x)]_{\alpha }\gamma y+S(x)\gamma [y,T(x)]_{\alpha },S(x)]_{\beta }
\\&=&S(x)\gamma [[S(x),T(x)]_{\alpha },y]_{\beta }+[[S(x),T(x)]_{\alpha },S(x)]_{\beta }\gamma y\\&&+[S(x),T(x)]_{\alpha }\gamma [y,S(x)]_{\beta }+T(x)\gamma [[S(x),y]_{\alpha },S(x)]_{\beta }\\&&+[T(x),S(x)]_{\beta }\gamma [S(x),y]_{\alpha }+[[S(x),T(x)]_{\alpha },S(x)]_{\beta }\gamma y\\&&+[S(x),T(x)]_{\alpha }\gamma [y,S(x)]_{\beta }+S(x)\gamma [[y,T(x)]_{\alpha },S(x)]_{\beta }\end{eqnarray*} 
Therefore, we have \begin{eqnarray}0&=&S(x)\gamma [[S(x),T(x)]_{\alpha },y]_{\beta }+3[S(x),T(x)]_{\alpha }\gamma [y,S(x)]_{\beta }\nonumber\\&&+T(x)\gamma [[S(x),y]_{\alpha },S(x)]_{\beta }+S(x)\gamma [[y,T(x)]_{\alpha },S(x)]_{\beta }\end{eqnarray} Replacing $y$ by $y\delta S(x)$ in the above relation, we have \begin{eqnarray*}0&=&S(x)\gamma [[S(x),T(x)]_{\alpha },y\delta S(x)]_{\beta }+3[S(x),T(x)]_{\alpha }\gamma [y\delta S(x),S(x)]_{\beta }\nonumber\\&&+T(x)\gamma [[S(x),y\delta S(x)]_{\alpha },S(x)]_{\beta }+S(x)\gamma [[y\delta S(x),T(x)]_{\alpha },S(x)]_{\beta }\\&=&S(x)\gamma [[S(x),T(x)]_{\alpha },y]_{\beta }\delta S(x)+S(x)\gamma y\delta [[S(x),T(x)]_{\alpha },S(x)]_{\beta }\\&&+3[S(x),T(x)]_{\alpha }\gamma [y,S(x)]_{\beta }\delta S(x)+T(x)\gamma [[S(x),y]_{\alpha }\delta S(x),S(x)]_{\beta }\\&&+S(x)\gamma [[y,T(x)]_{\alpha }\delta S(x)+y\delta [S(x),T(x)]_{\alpha },S(x)]_{\beta } 
\\&=&S(x)\gamma [[S(x),T(x)]_{\alpha },y]_{\beta }\delta S(x)+3[S(x),T(x)]_{\alpha }\gamma [y,S(x)]_{\beta }\delta S(x)\\&&+T(x)\gamma [[S(x),y]_{\alpha },S(x)]_{\beta }\delta S(x)+S(x)\gamma [[y,T(x)]_{\alpha },S(x)]_{\beta }\delta S(x)\\&&+S(x)\gamma [y,S(x)]_{\beta }\delta [S(x),T(x)]_{\alpha }+S(x)\gamma y\delta [[S(x),T(x)]_{\alpha },S(x)]_{\beta }\end{eqnarray*}  
Thus we have according to (36) and (38), \begin{eqnarray*}0&=&S(x)\gamma [y,S(x)]_{\beta }\delta [S(x),T(x)]_{\alpha }\end{eqnarray*} which can be written in the form  \begin{eqnarray}S(x)\gamma y\beta S(x)\delta [S(x),T(x)]_{\alpha }&=&S(x)\gamma S(x)\beta y\delta [S(x),T(x)]_{\alpha }\end{eqnarray} Putting in the above relation $y=T(x)\omega y$, we have \begin{eqnarray}S(x)\gamma T(x)\omega y\beta S(x)\delta [S(x),T(x)]_{\alpha }&=&S(x)\gamma S(x)\beta T(x)\omega y\delta [S(x),T(x)]_{\alpha }\end{eqnarray} On the other hand, left multiplication of (39) by $T(x)\omega $, we have \begin{eqnarray}\lefteqn{T(x)\omega S(x)\gamma y\beta S(x)\delta [S(x),T(x)]_{\alpha }=}\nonumber\\&&T(x)\omega S(x)\gamma S(x)\beta y\delta [S(x),T(x)]_{\alpha }\end{eqnarray}
Subtracting (41) from (40), we obtain \begin{eqnarray*}0&=&[S(x),T(x)]_{\gamma }\omega y\beta S(x)\delta [S(x),T(x)]_{\alpha }-[S(x)\gamma S(x),T(x)]_{\beta }\omega y\delta [S(x),T(x)]_{\alpha }\\&=&[S(x),T(x)]_{\gamma }\omega y\beta S(x)\delta [S(x),T(x)]_{\alpha }\\&&-([S(x),T(x)]_{\beta }\gamma S(x)+S(x)\gamma [S(x),T(x)]_{\beta }\omega y\delta [S(x),T(x)]_{\alpha }\end{eqnarray*}  According to the requirement of the theorem one can replace in the above calculation $[S(x),T(x)]_{\beta }\gamma S(x)$  by $S(x)\gamma [S(x),T(x)]_{\beta }$ which gives \begin{eqnarray*}\lefteqn{[S(x),T(x)]_{\gamma }\omega y\beta S(x)\delta [S(x),T(x)]_{\alpha }}\\&&=2S(x)\gamma [S(x),T(x)]_{\beta }\omega y\delta [S(x),T(x)]_{\alpha }\end{eqnarray*}
Left multiplication of the above relation by $S(x)\rho $ gives \begin{eqnarray}\lefteqn{S(x)\rho [S(x),T(x)]_{\gamma }\omega y\beta S(x)\delta [S(x),T(x)]_{\alpha }}\nonumber\\&&=2S(x)\rho S(x)\gamma [S(x),T(x)]_{\beta }\omega y\delta [S(x),T(x)]_{\alpha }\end{eqnarray} On the otherhand, putting $y=[S(x),T(x)]\rho y$ in (39), we have \begin{eqnarray}\lefteqn{S(x)\gamma [S(x),T(x)]_{\omega }\rho y\beta S(x)\delta [S(x),T(x)]_{\alpha }}\nonumber\\&&=S(x)\gamma S(x)\beta [S(x),T(x)]_{\omega }\rho y\delta [S(x),T(x)]_{\alpha }\end{eqnarray} Combining (42) with (43), we obtain \begin{eqnarray*}0&=&S(x)\rho [S(x),T(x)]_{\gamma }\omega y\beta S(x)\delta [S(x),T(x)]_{\alpha }\end{eqnarray*} Hence by semiprimeness of $M$, we obtain \begin{eqnarray}S(x)\delta [S(x),T(x)]_{\alpha }&=&0\end{eqnarray} 
From (44) and the assumption of the theorem, we have \begin{eqnarray*}[S(x),T(x)]_{\alpha }\delta S(x)&=&0\end{eqnarray*} The rest of the proof goes through in the same way as in the proof of the Theorem-2.1.

\end{document}